# Characterizations and Infinite Divisibility of Extended COM-Poisson Distribution

**Huiming Zhang**

School of Mathematics and Statistics, Central China Normal University, Wuhan, China

**Email address:** zhanghuiming@mails.ccnu.edu.cn

**Abstract:** This article provides some characterizations of extended COM-Poisson distribution: conditional distribution given the sum, functional operator characterization (Stein identity). We also give some conditions such that the extended COM-Poisson distribution is infinitely divisible, hence some subclass of extended COM-Poisson distributions are discrete compound Poisson distribution.

## 1. Extended COM-Poisson Distribution

The Conway-Maxwell-Poisson distribution (COM-Poisson distribution) was firstly briefly introduced by Conway and Maxwell (1962) for modeling of queuing systems with state-dependent service times, see Shmueli et al. (2005) for details. The probability mass function (p.m.f.) is given by

$$P(X=k) = \frac{\lambda^k}{(k!)^\nu} \cdot \frac{1}{Z(\lambda,\nu)}, (k=0,1,2,\cdots),$$

where $\lambda, \nu > 0$ and $Z(\lambda,\nu) = \sum_{i=0}^{\infty} \frac{\lambda^i}{(i!)^\nu}$.

Chakraborty(2015) introduced the extended COM-Poisson distribution from Imoto (2014), Chakraborty and Ong (2014):

**Definition A**: *A r.v. X is said to follow the extended COM-Poisson distribution with parameters $(\nu, p, \alpha, \beta)$ [ $X \sim$ ECOMP $(\nu, p, \alpha, \beta)$ ] if its p.m.f. is given by*

$$P(X=k) = \frac{\frac{[\Gamma(\nu+k)]^\beta}{[\Gamma(\nu)]^\beta (k!)^\alpha} p^k}{\sum_{i=0}^{\infty} \frac{[\Gamma(\nu+i)]^\beta}{[\Gamma(\nu)]^\beta (i!)^\alpha} p^i} =: \frac{\frac{[\Gamma(\nu+k)]^\beta}{[\Gamma(\nu)]^\beta (k!)^\alpha} p^k}{{}_1S^\beta_{\alpha-1}(\nu;1;p)}, (k=0,1,\cdots)$$

where the parameter space is

$$(\nu \geq 0, p > 0, \alpha > \beta \geq 0) \cup (\nu > 0, 0 < p < 1, \alpha = \beta \geq 0).$$

The extended COM-Poisson includes many COM type distributions. For example, COM-Poisson, COM-negative binomial etc., see Chakraborty(2015). It is easy to obtain the p.m.f. of recurrence relation of the ECOMP $(\nu, p, \alpha, \beta)$

$$\frac{P(X=k+1)}{P(X=k)} = \frac{p(\nu+k)^\beta}{(k+1)^\alpha} \qquad (1)$$

with $P(X=0) = [{}_1S^\beta_{\alpha-1}(\nu;1;p)]^{-1}$.

The recurrence relation (1) will be useful for arriving some characterizations and properties of the extended COM-Poisson distribution.

There are many statistical models pertain to COM-Poisson for modeling many types of count data, such as COM-Poisson regression (see Guikema and Goffelt (2008)), COM-Poisson in survival analysis (see Rodrigues et al. (2009), Sellers and Shmueli (2010)), COM-Poisson INGARCH time series models (see Zhu (2012)), COM-Poisson distribution chart in statistical process control (see Saghir et al. (2009)), zero-inflated COM-Poisson distribution (see Barriga and Phillips (2014)), fitting data from actuarial science (see Khan and Khan (2010)).

Undoubtedly, extended COM-Poisson distribution with 4 parameters are more flexible than COM-Poisson distribution, negative binomial distribution and so on. It contains a wide range of statistical properties (log-concave and log-convex, under- and over dispersed). Especially, for the phenomenon of over-/equi-/underdispersion, there are many applications in applied statistics, including public health, medicine, and epidemiology, see the literatures review by Kokonendji (2014). So we can consider statistical models above in terms of extended COM-Poisson distribution in the future.

### 1.1. Functional Operator Characterization

We know the p.m.f. of ECOMP $(\nu, p, \alpha, \beta)$ has a simple recurrence relation given by (1). Brown and Xia (2001) studied a very large class of stationary distribution of birth-death process with arrival rate $\lambda_k$ and service rate $\mu_k$ by the recursive formula:

$$\frac{P(X=k+1)}{P(X=k)} = \frac{\lambda_k}{\mu_{k+1}}, (k=0,1,\cdots).$$

Thus we can construct that arrival rate

$$\lambda_k = (\nu+k)^\beta \lambda, (k=1,2,\cdots)$$

and service rate

$$\mu_k = k^\alpha \mu, (k=1,2,\cdots)$$

to characterize the ECOMP $(\nu, p, \alpha, \beta)$, $p = \lambda/\mu$, see also Daly and Gaunt (2015). This will be useful for functional characterization in Stein's method. The following theorem presents a characterization for the ECOMP $(\nu, p, \alpha, \beta)$ distributions.

**Theorem A**: (Stein identity) *Let* N={0,1,2,$\cdots$}, *the r.v. X has distribution* ECOMP $(\nu, p, \alpha, \beta)$ *iff the equation*

$$E[X^\alpha g(X) - p(X+\nu)^\beta g(X+1)] = 0$$

*holds for any bounded function* $g: N \to R$.

The Theorem A can be directly obtained from the work of Brown and Xia (2001). On the one hand, we just check the expectation is 0, and then the sufficiency is true. On the other hand, to show the necessary, just take function $g$ to be the indicator of $\{k\}$, hence we have the recurrence relation (1).

### *1.2. Conditional Distribution Characterization*

For the two independent r.v.'s

$X \sim$ ECOMP $(\nu_1, p, \alpha, \beta)$ and $Y \sim$ ECOMP $(\nu_2, p, \alpha, \beta)$, the distribution of sum $S = X + Y = s$ is

$$P(X+Y=s) = \sum_{k=0}^{s} P(X=k)P(Y=s-k)$$
$$= \frac{[\Gamma(\nu_1)\Gamma(\nu_2)]^{-\beta}}{{}_1S_{\alpha-1}^\beta(\nu_1;1;p)\,{}_1S_{\alpha-1}^\beta(\nu_2;1;p)} \sum_{k=0}^{s} \frac{\{\Gamma(\nu_1+k)\Gamma(\nu_2+s-k)\}^\beta}{[k!(s-k)!]^\alpha} p^s$$

So, the conditional distribution of $X=k \mid X+Y=s$ is

$$P(X=k|S=s) = \frac{\{\Gamma(\nu_1+k)\Gamma(\nu_2+s-k)\}^\beta}{[k!(s-k)!]^\alpha} \Big/ \sum_{x=0}^{s} \frac{\{\Gamma(\nu_1+x)\Gamma(\nu_2+s-x)\}^\beta}{[x!(s-x)!]^\alpha}.$$

We naturally define the extended negative hypergeometric distribution with p.m.f.

$$P(Z=k) = \frac{\{\Gamma(\nu_1+k)\Gamma(\nu_2+s-k)\}^\beta}{[k!(s-k)!]^\alpha} \Big/ E(s, \nu_1, \nu_2, \alpha, \beta), \quad (k=0,1,\cdots,s) \quad (2)$$

where $E(s, \nu_1, \nu_2, \alpha, \beta)$ is the normalization constant.

And we denote (2) as $Z \sim$ ENHG$(s, \nu_1, \nu_2, \alpha, \beta)$.

Next, we cite a general result by Patil and Seshadri (1964) for characterizing a large class of discrete distributions, see also Kagan et al. (1973). And we obtain the conditional distribution characterization of ECOMP $(\nu_1, p, \alpha, \beta)$.

**Lemma A**: *Let X and Y be independent discrete r.v.'s and two dimensional function be*

$$c(x, x+y) = P(X=x|X+Y=x+y).$$

*If*

$$\frac{c(x+y, x+y)c(0, y)}{c(x, x+y)c(y, y)} = \frac{h(x+y)}{h(x)h(y)}$$

*where h is a nonnegative function, then*

$$f(x) = f(0)h(x)e^{ax}, g(y) = g(0)\frac{h(y)h(0,y)}{c(y,y)} e^{ay}$$

*where*

$$0 < f(x) = P(X=x), 0 < g(y) = P(Y=y).$$

**Theorem B**: *Let X and Y be the independent discrete r.v. with*

$$P(X=x) = f(x) > 0 \text{ and } P(Y=y) = g(y) > 0.$$

*If the conditional distribution of $X = x \mid X+Y=s$ is extended negative hypergeometric distribution* ENHG$(s, \nu_1, \nu_2, \alpha, \beta)$, *then*

$X \sim$ ECOMP $(\nu_1, p, \alpha, \beta)$ *and* $Y \sim$ ECOMP $(\nu_2, p, \alpha, \beta)$.

Proof: By using Lemma A and the definition of ENHG$(s, \nu_1, \nu_2, \alpha, \beta)$, we have

$$c(x, x+y) = \frac{\{\Gamma(\nu_1+x)\Gamma(\nu_2+x+y-x)\}^\beta}{[x!(x+y-x)!]^\alpha} \Big/ E(s, \nu_1, \nu_2, \alpha, \beta).$$

Then

$$\frac{c(x+y, x+y)c(0,y)}{c(x, x+y)c(y,y)} = \frac{\frac{\{\Gamma(\nu_1+x+y)\Gamma(\nu_2)\}^\beta}{[(x+y)!]^\alpha} \frac{\{\Gamma(\nu_1)\Gamma(\nu_2+y)\}^\beta}{(y!)^\alpha}}{\frac{\{\Gamma(\nu_1+x)\Gamma(\nu_2+x+y-x)\}^\beta}{[x!(x+y-x)!]^\alpha} \frac{\{\Gamma(\nu_1+y)\Gamma(\nu_2)\}^\beta}{(y!)^\alpha}}$$
$$= \frac{[\Gamma(\nu_1+x+y)]^\beta}{[\Gamma(\nu_1)]^\beta[(x+y)!]^\alpha} \Big/ \frac{[\Gamma(\nu_1+x)]^\beta}{[\Gamma(\nu_1)]^\beta(x!)^\alpha} \frac{[\Gamma(\nu_1+y)]^\beta}{[\Gamma(\nu_1)]^\beta(y!)^\alpha}$$

So, we have

$$h(x) = \frac{[\Gamma(\nu_1+x)]^\beta}{[\Gamma(\nu_1)]^\beta(x!)^\alpha}$$

and

$$\frac{h(y)c(0,y)}{c(y,y)} = \frac{[\Gamma(\nu_1+y)]^\beta}{[\Gamma(\nu_1)]^\beta(y!)^\alpha} \frac{\{\Gamma(\nu_1)\Gamma(\nu_2+y)\}^\beta}{\{\Gamma(\nu_1+y)\Gamma(\nu_2)\}^\beta} = \frac{[\Gamma(\nu_2+y)]^\beta}{[\Gamma(\nu_2)]^\beta(y!)^\alpha}.$$

Let $p = e^a$, then

$$P(X=x) = \frac{[\Gamma(\nu_1+x)]^\beta}{[\Gamma(\nu_1)]^\beta(x!)^\alpha} p^x P(X=0),$$

$$P(Y=y) = \frac{[\Gamma(\nu_2+y)]^\beta}{[\Gamma(\nu_2)]^\beta(y!)^\alpha} p^y P(Y=0).$$

Hence we get the result. □

## 2. Log-Convex and Infinitely Divisible

For a discrete r.v. $X$ with p.g.f.

$$P_X(z) = \sum_{k=0}^{\infty} P(X=k)z^k, (|z| \leq 1)$$

We say that $X$ is infinitely divisible if

$$P_X(z) = [P_{X_n}(z)]^n, (\forall n \in N, |z| \leq 1)$$

where $P_{X_n}(z)$ is the p.g.f. of a certain discrete r.v. $X_n$.

Feller's characterization of the infinite divisible discrete r.v. (see Feller (1971)) shows that a discrete r.v. is infinitely divisible iff its distribution is a discrete compound Poisson distribution with following p.g.f.

$$G(z) = \sum_{k=0}^{\infty} p_k z^k = e^{\sum_{i=1}^{\infty} \alpha_i \lambda (z^i - 1)}, (|z| \leq 1)$$

where $\alpha_i \lambda$ is the Lévy measure(or parameters) of a infinitely divisible distribution. It satisfies

$$\sum_{i=1}^{\infty} \alpha_i = 1, \alpha_i \geq 0, \lambda > 0.$$

Actually, the discrete compound Poisson distribution is the probability distribution of the sum of a number of iid non-negative integer-valued r.v.'s, where the sum is Poisson-distributed r.v..

On the one hand, the discrete compound Poisson(DCP) distributed r.v. can be represented as the sum of $n$ i.i.d. r.v.'s

$$X = Y_1 + Y_2 + \cdots + Y_N,$$

where $N$ and $Y_i (i=0,1,2,\cdots)$ are independent with

$$P(Y_1 = i) = \alpha_i, P(N = n) = \frac{\lambda^n}{n!} e^{-\lambda}.$$

On the other hand, $X$ can be decomposed as sum of weighted Poisson:

$$X = Z_1 + 2Z_2 + \cdots + iZ_i + \cdots,$$

where $Z_i (i=1,2,\cdots)$ are independently Poisson distributed with parameter $\alpha_i \lambda$.

For the detailed theoretical and applied treatment of discrete infinitely divisible and discrete compound Poisson, we refer readers to section 2 of Steutel and van Harn (2003), section 9.3 of Johnson et al. (2005), Zhang et al. (2014), Zhang et al. (2013).

A discrete r.v. $X$ with $p_k = P(X=k)$ has log-concave (log-convex) p.m.f. if

$$\frac{p_{k+1} p_{k-1}}{p_k^2} = \frac{p_{k+1}}{p_k} \Big/ \frac{p_k}{p_{k-1}} \leq (\geq) 1, \quad (k \geq 1).$$

Steutel (1970) proved that all log-convex discrete distributions are infinitely divisible, see also Steutel and van Harn (2003). So, it is easy to get infinite divisibility of ECOMP $(v, p, \alpha, \beta)$ when $(v, \alpha, \beta)$ satisfy some conditions.

**Theorem C**: *The ECOMP $(v, p, \alpha, \beta)$ infinitely divisible distribution when $(v, \alpha, \beta)$ satisfy the following conditions:*

$$\frac{1}{2^\alpha}\left(1 + \frac{1}{v}\right)^\beta \geq 1.$$

Proof: From (1), we have

$$\frac{P(X = k+1)}{P(X = k)} \Big/ \frac{P(X = k)}{P(X = k-1)} = \frac{k^{\alpha-\beta}}{(k+1)^{\alpha-\beta}} \frac{(k^2 + vk)^\beta}{(k^2 + vk + v - 1)^\beta}, (k=1,2,\cdots).$$

Then, the equality

$$\frac{P(X = k+1)}{P(X = k)} \Big/ \frac{P(X = k)}{P(X = k-1)} \geq 1$$

holds if

$$\frac{k^{\alpha-\beta}}{(k+1)^{\alpha-\beta}} \frac{(k^2 + vk)^\beta}{(k^2 + vk + v - 1)^\beta} \geq \frac{1}{2^{\alpha-\beta}} \frac{(k^2 + vk)^\beta}{(k^2 + vk + v - 1)^\beta} \geq 1 \text{ for } k=1,2,\cdots.$$

Notice that $\frac{(k^2 + vk)^\beta}{(k^2 + vk + v - 1)^\beta}$ goes to its minimum $\left(\frac{1+v}{2v}\right)^\beta$ as $k = 1$. Hence, we need

$$\frac{1}{2^\alpha}\left(1 + \frac{1}{v}\right)^\beta \geq 1$$

to make sure the infinite divisibility of ECOMP $(v, p, \alpha, \beta)$. □

Applying the recurrence relation (Lévy-Adelson-Panjer recursion) of p.m.f. of DCP distribution,

$$P_{n+1} = \frac{\lambda}{n+1}[\alpha_1 P_n + 2\alpha_2 P_{n-1} + \cdots + (n+1)\alpha_{n+1} P_0], (P_0 = e^{-\lambda}, n=0,1,\cdots),$$

see Buchmann (2003) and the Remark 1 in Zhang et al. (2014), therefore DCP case of ECOMP $(v, p, \alpha, \beta)$ has the alternative recurrence relation.

The parameters $\alpha_i (i=1,2,\cdots)$ and $\lambda$ of DCP case of ECOMP $(v, p, \alpha, \beta)$ are determined by the following systems of equation:

$$\frac{\{\Gamma(v+n+1)\}^\beta p^{n+1}}{[(n+1)!]^\alpha} = \frac{\lambda}{n+1}[\alpha_1 \frac{\{\Gamma(v+n)\}^\beta p^n}{(n!)^\alpha} +$$
$$2\alpha_2 \frac{\{\Gamma(v+n-1)\}^\beta p^{n-1}}{[(n-1)!]^\alpha} + \cdots + (n+1)\alpha_{n+1}\{\Gamma(v)\}^\beta], (n=0,1,\cdots)$$

where $\lambda = \ln {}_1 S_{\alpha-1}^\beta (v; 1; p)$ .

The is ECOMP $(v, p, \alpha, \beta)$ log-concave distribution if $v, \alpha, \beta$ satisfy the inequality

$$\frac{P(X = k+1)}{P(X = k)} \Big/ \frac{P(X = k)}{P(X = k-1)} = \frac{k^{\alpha-\beta}}{(k+1)^{\alpha-\beta}} \frac{(k^2 + vk)^\beta}{(k^2 + vk + v - 1)^\beta} \leq 1 \text{ for } k=1,2,\cdots.$$

Note that $v \geq 1$ ensures validity of the inequality above. For more theoretical results about general log-concavity of discrete distributions, see Balabdaoui et al. (2013), Saumard and Wellner (2014).

## Acknowledgments

In the author's future articles, more new characterizations and properties related to this notes in aspects of COM-Poisson distribution will appear, and Extended COM-Poisson distribution in aspects of related statistical models are in progress. If you are interested it, please contact me without hesitation.